\documentclass[12pt,reqno]{amsart}

\newtheorem{theorem}{Theorem}

\newtheorem{proposition}[theorem]{Proposition}
\newtheorem{corollary}[theorem]{Corollary}

\theoremstyle{definition}

\newtheorem{example}{Example}
\newtheorem{claim}{Claim}

\theoremstyle{remark}

\newcommand{\ltwo}{L^2({\mathbb R})}
\newcommand{\bbz}{\mathbb{Z}}

\newcommand{\Cal}{\mathcal}
\newcommand{\Bb}{\mathbb}

\newcommand{\wo}{\widehat{W}_0}

\begin{document}

\title[Translation Invariant Wavelets]{Wavelets having the Translation Invariance Property of Order n}
\author{Sharon Schaffer}
\address{Department of Mathematics, University of Colorado,
Boulder, CO, 80309-0395} 
\email{schaffs@euclid.colorado.edu}
\author{Eric Weber}
\address{Department of Mathematics, Texas A\&M University, College Station, TX 77843-3368}
\email{weber@math.tamu.edu}
\begin{abstract}
All wavelets can be associated to a multiresolution like structure, i.e.~an increasing sequence of subspaces of $\ltwo$.  We consider the interaction of a wavelet and the translation operator in terms of which of the subspaces in this multiresolution like structure are invariant under the translation operator.  This action defines the notion of the translation invariance property of order $n$.  In this paper we show that wavelets of all levels of translation invariance exist, first for the classic case of dilation by 2, and then for arbitrary integral dilation factors.
\end{abstract}
\maketitle

\section{Introduction}

The study of the interplay between wavelets and operator theory has proven useful recently.  This paper considers part of that study, in particular the interaction of certain wavelets and the translation operator.  Every wavelet can be associated to a GMRA, which is a nested sequence of subspaces of $\ltwo$ with other properties.  The central subspace, called the \emph{core subspace} has the property that it is, along with the larger spaces in this ``ladder'', invariant under integer translations.  In this paper, we shall show the existence of wavelets such that some, but not all, of the smaller spaces in this ladder are also invariant under integer translations.  In this context, these wavelets can be considered generalizations of MSF wavelets.  We begin with several definitions.

A wavelet for dilation $d$ is a square integrable function such that the collection
\[ \{d^{\frac{n}{2}} \psi(d^n x - l): n,\ l \in \Bb{Z} \} \]
forms an orthonormal basis of $\ltwo$.  Define unitary operators on $\ltwo$ by $Tf(x) = f(x -1)$ and $Df(x) = \sqrt{d}f(dx)$ for $d \in \mathbb{Z}, \: d \geq 2$, called the translation and dilation operator, respectively.  As in \cite{BMM}, define a nonsingular \emph{Generalized Multiresolution Analysis} (GMRA) as a sequence of subspaces $\{V_j\}_j$ with the following five properties:
\begin{enumerate}
\item $V_j \subset V_{j+1}$,
\item $\bigcup_{j \in \bbz} V_j$ is dense in $\ltwo$,
\item $\bigcap_{j \in \bbz} V_j = \{0\}$,
\item $D(V_{j}) = V_{j+1}$
\item $V_0$ is invariant under $T^l$ for $l \in \Bb{Z}$.
\end{enumerate}
Given a wavelet $\psi$, define $V_j$ as the closed span of $\{D^n T^l \psi: n < j, \ l \in \Bb{Z}\}$.  It is routine to verify that this sequence of subspaces of $\ltwo$ form a GMRA.  Lastly, define a dual sequence of closed subspaces, denoted $W_j$, by $V_{j+1}=V_j \oplus W_j$.  Note that $W_0 = \overline{span}\{T^l \psi: l \in \Bb{Z}\}$, i.e.~$W_0$ is where $\psi$ lives.

The commutation relation between $D$ and $T$, $TD = DT^d$, yields that for any $j$, if $V_j$ is invariant under $T$, then so is $V_{j+1}$.  Hence, for a GMRA, the subspaces $V_{j}$ are invariant under $T$ for $j \geq 0$.  We will show by construction that, given a negative integer $n$, there exists a wavelet such that its associated GMRA has the property that the subspaces $V_j$ are invariant under $T^l$ if and only if $j \geq n$.  Such a wavelet is said to have the translation invariance property of order $n$.

We shall denote by $\Cal{L}_n$ the collection of all wavelets such that $V_{-n}$ is invariant under integral translations.  Note that these collections are nested:
\[
\Cal{L}_0 \supset \Cal{L}_1 \supset \Cal{L}_2 \ldots \supset \Cal{L}_n \supset \Cal{L}_{n+1} \supset \ldots \supset \Cal{L}_{\infty}
\]
Define $\Cal{M}_n = \Cal{L}_n \setminus \Cal{L}_{n+1}$, with $\Cal{M}_{\infty} = \Cal{L}_{\infty}$, the collection of all wavelets such that $V_{j}$ is invariant under translations for all $j \in \Bb{Z}$.  It is proven in \cite{W2} that the Minimally Supported Frequency (MSF) wavelets coincide with $\Cal{M}_{\infty}$.  A MSF wavelet is a wavelet $\psi$ such that $|\hat{\psi}|$ is the indicator function of some set, usually called a wavelet set.  MSF wavelets shall be a starting point for constructing the wavelets we seek.

The goal of this paper is to prove the following theorem:
\begin{theorem} \label{T:Main}
For all positive integers $n$, there exists a wavelet with translation invariance property of order $n$, i.e. the collections $\Cal{M}_n$ are non-empty.  This is true for all $n \in \Bb{N}$ and for all positive integer dilation factors.
\end{theorem}
The proof is the content of section~\ref{S:Proof}.

For the purposes of this paper, we shall say that a set $E \subset \Bb{R}$ is \emph{partially self-similar} with respect to $\alpha \in \Bb{R}$ if there exists a set $F$ of non-zero measure such that both $F$ and $F + \alpha$ are subsets of $E$.  Additionally, if $G,H$ are two subsets of $\Bb{R}$, we shall say that $G$ is $2 \pi$ translation congruent to $H$ if there exists a measurable partition $G_n$ of $G$ such that the collection $\{G_n + 2 n \pi: n \in \Bb{Z}\}$ forms a partition of $H$, modulo sets of measure zero.  Define a mapping $\tau: \Bb{R} \to [0, 2\pi)$ such that $\tau(x) - x = 2 \pi k$ for some integer $k$.

\section{A Characterization}

In this section we shall develop a characterization of wavelets that have the translation invariance property of order $n$.  This characterization comes from considering the support of the Fourier transform of the wavelet in question.  Fix a dilation factor $d$.  Denote by $T_{\alpha}$ the unitary operator $T_{\alpha}f(x) = f(x - \alpha)$.  $T$ is to be understood as $T_1$.  Note that $\widehat{T_{\alpha}} = M_{e^{-i \alpha \xi}}$.  We shall consider the groups of translations $\Cal{G}_n = \{ T_{\frac{m}{d^n}}: m \in \Bb{Z} \}$.  The following proposition will allow us to consider the action of translation operators on $W_0$ instead of on the $V_j$'s.

\begin{proposition} \label{L:VG}
The following are equivalent:
\begin{enumerate}
\item the space $V_{-n}$ for $\psi$ is invariant under the action of translations,
\item $V_0$ is invariant under the action of $T_{\frac{m}{d^n}}$,
\item $W_0$ is invariant under the action of $T_{\frac{m}{d^n}}$.
\end{enumerate}
\end{proposition}
\begin{proof}
By definition, $f \in V_{-n}$ if and only if $D^n f \in V_0$; this fact combined with the commutation relation $DT = T_{\frac{1}{d}}D$ establishes the equivalence of $1$ and $2$.  Furthermore, these two facts show that if $V_0$ is invariant under the action of $T_{\frac{m}{d^n}}$, then so is $V_1$, whence $W_0$ is as well.

Finally, suppose that $W_0$ is invariant under $T_{\frac{m}{d^n}}$.  Then for $k > 0$, $W_k$ is invariant under $T_{\frac{m}{d^n}}$, whence $V_0^{\perp} = \oplus_{k=0}^{\infty} W_k$ is invariant under $\Cal{G}_n$.  It follows that $V_0$ is invariant under $\Cal{G}_n$.
\end{proof}

If $f \in W_0$, then we can write $f = \sum_{k \in \Bb{Z}} c_k T^k \psi$, so taking the Fourier transform of both sides yields $\hat{f} = g \hat{\psi}$ for some $g \in L^2([0, 2\pi))$.  Hence, we can describe $\wo$ by $\{ g(\xi) \hat{\psi}(\xi) : g \in L^2([0, 2\pi)) \}$.  Suppose $\psi$ is a wavelet that is in $\Cal{L}_n$, by proposition~\ref{L:VG}, this is equivalent to $T_{\frac{m}{d^n}} f \in W_0$.  By taking the Fourier Transform, we have $e^{-i \frac{m}{d^n} \cdot} \hat{\psi} \in \wo$.  Hence, $\psi \in \Cal{L}_n$ is equivalent to the condition that $e^{-i \frac{m}{d^n} \xi} \hat{\psi}(\xi) = g(\xi) \hat{\psi}(\xi)$.  For ease of notation, define:
\[\Cal{E}_{\psi} = \{k \in \Bb{Z}:\ supp(\hat{\psi}) \text{ is partially self similar w.r.t. } 2\pi k \}. \]

\begin{theorem}
Let $\psi$ be a wavelet.  Then $\psi \in \Cal{L}_n$ if and only if every element of $\Cal{E}_{\psi}$ is divisible by $d^n$.
\end{theorem}
\begin{proof}[Only If]
We prove by contradiction.  Let $\psi \in \Cal{L}_n$.  Hence,
\[ e^{-i \frac{1}{d^n} \xi} \hat{\psi}(\xi) = g(\xi) \hat{\psi}(\xi) \]
for some $g \in L^2([0, 2\pi))$.  Let $E = supp(\hat{\psi})$, for $\xi \in E$, $e^{-i \frac{1}{d^n} \xi} = g(\xi)$.  Suppose that $k \in \Cal{E}_{\psi}$ is not divisible by $d^n$, and $F$ is a set of non-zero measure such that both $F$ and $F+2k\pi$ are subsets of $E$.  Let $\xi \in F$; we have,
\begin{align*}
e^{-i \frac{\xi +2k\pi}{d^n}} \hat{\psi}(\xi + 2 k\pi) &= g(\xi+ 2 k\pi) \ \hat{\psi}(\xi + 2k\pi) \\
e^{-i \frac{\xi}{d^n}} e^{-i \frac{2k\pi}{d^n}} \hat{\psi}(\xi + 2 k\pi) &=  g(\xi) \hat{\psi}(\xi + 2k\pi) \\
e^{-i \frac{2k\pi}{d^n}} &= 1
\end{align*}
a contradiction.

\noindent \emph{If.}
It suffices to show that 
\[ e^{-i \frac{1}{d^n} \xi} \hat{\psi}(\xi) = g(\xi) \hat{\psi}(\xi) \] for some $g \in L^2([0, 2\pi))$.  Again, let $E = supp(\hat{\psi})$.

Let $F \subset E$ be such that $\tau:F \to [0, 2\pi)$ is a bijection.  (It can be easily shown that $\tau:E \to [0, 2\pi)$ is a surjection.)  The injective property of $\tau$ can be assured in the following manner: for each $\xi \in [0, 2\pi)$, define the set $Z_{\xi} = \{m_{\xi} \in \bbz: \xi + 2 m_{\xi} \pi \in E\}$, then for $\xi$ choose $k_{\xi}$ to be 0 if $\xi \in E$, if not, choose $k_{\xi} = min\{m > 0: m \in Z_{\xi}\}$, else choose $k_{\xi} = max\{m < 0: m \in Z_{\xi}\}$.  Let $F = \{\xi + 2 k_{\xi} \pi : \xi \in [0, 2\pi) \}$.  Note that by construction, $F$ is $2\pi$ translation congruent to $[0, 2\pi)$.  Hence,
\[ e^{-i \frac{1}{d^n} \xi} \chi_{F}(\xi) = g(\xi) \]
where $g(\xi) \in L^2(F)$ and is $2 \pi$ periodic.  Thus, for $\xi \in F$, 
\[ e^{-i \frac{1}{d^n} \xi} \hat{\psi}(\xi) = g(\xi) \hat{\psi}(\xi). \]

For almost any $\xi \in E \setminus F$, there exists a $\xi' \in F$ and an integer $l_{\xi}$ such that $\xi - \xi' = 2 l_{\xi} \pi$.  Moreover, by hypothesis, $l_{\xi}$ is a multiple of $d^n$, since $l_{\xi} \in \Cal{E}_{\psi}$.  Since $e^{-i \frac{1}{d^n} \xi}$ is $d^n \pi$ periodic, we have that for $\xi \in E \setminus F$,
\begin{align*}
e^{-i \frac{1}{d^n} \xi} \hat{\psi}(\xi) &=  e^{-i \frac{1}{d^n} (\xi' + 2l_{\xi}\pi)} \hat{\psi}(\xi' + 2l_{\xi}\pi) \\
	&= e^{-i \frac{1}{d^n} \xi'} \hat{\psi}(\xi' + 2l_{\xi}\pi) \\
	&= g(\xi') \hat{\psi}(\xi' + 2l_{\xi}\pi) \\
	&= g(\xi) \hat{\psi}(\xi).
\end{align*}
This completes the proof.
\end{proof}

As an immediate corollary, we have the following characterization of wavelets in $\Cal{M}_n$.
\begin{corollary} \label{C:1}
Let $\psi$ be a wavelet.  Then $\psi \in \Cal{M}_n$ if and only if every element of $\Cal{E}_{\psi}$ is divisible by $d^n$, but there exists an element of $\Cal{E}_{\psi}$ that is not divisible by $d^{n+1}$
\end{corollary}

\section{Wavelets of Higher Order Translation Invariance} \label{S:Proof}

In this section, we arrive at the crux of the paper, constructing wavelets with the translation invariance property of order $n$ .  When finished, we will have proven theorem~\ref{T:Main}.  The basic procedure is to consider a special class of \emph{wavelet multiplicity functions} (see \cite{BMM}), from which we can construct two generalized scaling sets.  Those in turn give rise to two wavelet sets which admit operator interpolation.  The resulting interpolated wavelet has the desired support in the frequency domain (proposition~\ref{L:VG}).  The multiplicity functions depend on whether the dilation factor is 2 or another integer; we shall first consider dilation by 2.

Recall that $\psi$ is a MSF wavelet if $|\hat{\psi}| = \chi_W$, where $W$ is referred to as a \emph{wavelet set}.  Such a wavelet shall be denoted by $\psi_W$.  We refer the reader to \cite{DL, FW, HW} for background information on wavelet sets.  A set $E$ is called a \emph{generalized scaling set} if the set $W = d E \setminus E$ is a wavelet set (\cite{BRS}).

We refer the reader to \cite[ch. 4]{DL} for a detailed description of operator interpolation.  We shall give here the necessary results concerning operator interpolation.  Let $\psi_{W_1}$ and $\psi_{W_2}$ be MSF wavelets.  By lemma 4.3 in \cite[pg. 41]{DL}, $W_1$ is $2 \pi$ translation congruent to $W_2$.  If $\sigma: W_1 \to W_2$ is determined by this translation congruence, then $\sigma$ can be extended to a measurable bijection of $\Bb{R}$ by defining $\sigma(x) = d^{-n} \sigma(d^n x)$ where $n$ is such that $d^n x \in W_1$.

If $\sigma$ is \emph{involutive}, i.e. $\sigma^2$ is the identity, and if $h_1$ and $h_2$ are measurable, essentially bounded, d-dilation periodic functions (i.e. $h_1(dx) = h_1(x)$), then $\psi$ defined by 
\[ \hat{\psi} = h_1 \hat{\psi}_{W_1} + h_2 \hat{\psi}_{W_2} \]
is again a wavelet provided the matrix
\begin{equation} \label{E:M}
\left( \begin{matrix}
h_1 & h_2 \\ h_2 \circ \sigma^{-1} & h_1 \circ \sigma^{-1} \end{matrix}
\right) \end{equation}
is unitary almost everywhere.  (Since $\sigma^{-1}$ is d-homogeneous, and the $h_i$'s are d-dilation periodic, in general it suffices to check this condition on $W_1$.)  Note that the interpolated wavelet $\psi$ has the property that $supp(\hat{\psi}) \subset W_1 \cup W_2$.  Furthermore, since $\sigma$ on $W_1$ is given by translations by integral multiples of $2 \pi$, $\sigma$ completely describes the partial self similarity of $W_1 \cup W_2$ with respect to multiples of $2 \pi$.

In order to find the wavelet sets and do the interpolation between them, we need the following theorem, which provides a technique for constructing generalized scaling sets.  The proof can be found in \cite{BM}.

\begin{theorem}[Merrill]  Given a multiplicity function $m:[-\pi, \pi) \to \{0, 1, \ldots, \infty\}$ which satisfies the consistency equation:
\begin{equation} \label{Eq:M}
m(\omega) + 1 = \sum_{i = 0}^{d - 1} m\left( \frac{\omega}{d} + \frac{2 \pi i}{d}\right)
\end{equation}
$m$ is associated with a (wavelet set $W$) wavelet if and only if $\exists$ a set E such that:
\begin{enumerate}
\item $m(\omega) = \sum_{k \in \bbz} \chi_E(\omega + 2 \pi k)$
\item $E \subseteq d E$
\item $\cup_{j \in \bbz} d^j E \textrm{ contains a neighborhood of 0.}$
\end{enumerate}
Indeed, $ W = d E \setminus E $ is such a wavelet set.
\end{theorem}

Baggett constructed an entire family of multiplicity functions for $d = 2$ that satisfy the consistency equation and produce wavelet sets, and hence wavelets.  Merrill discovered that these functions could be constructed using the idea of periodic points, and the breaks in these multiplicity functions occurred at these periodic points.  These periodic points are found by solving the following equation:
\[ d^k a = a + 2 \pi j\] where $ k - 1$ is the highest value of $m(\omega)$.

Solving for $a$ yields \[ a = \frac {2 \pi j}{d^k - 1}.\]  For a dilation factor of 2, Baggett constructed the following family of multiplicity functions, where $k \geq 2$ and the highest value of $m(\omega)$ is $k - 1$.

\begin{equation*} 
m_k^{2}(\omega) = \
\begin{cases}
1, &-\pi \leq \omega < \frac{-(2^k - 2)\pi}{2^k - 1};\\
0,  &\frac{-(2^k - 2)\pi}{2^k - 1} \leq \omega < \frac{-2^{k-1}\pi}{2^k - 1}; \\
k - j, &\frac{-2^{j}\pi}{2^k - 1} \leq \omega < \frac{- 2^{j-1} \pi}{2^k - 1};\\
k - 1, &\frac{-2\pi}{2^k - 1} \leq \omega < \frac{2 \pi}{2^k - 1};\\
k - j,  &\frac{2^{j-1} \pi}{2^k - 1} \leq \omega < \frac{2^j \pi}{2^k - 1};\\
0, &\frac{2^{k - 1}\pi}{2^k - 1} \leq \omega < \frac{(2^k - 2) \pi}{2^k - 1};\\
1, &\frac{(2^k - 2) \pi}{2^k - 1} \leq \omega < \pi,
\end{cases}
\end{equation*}
for $j = 2, 3, \ldots, k - 1$.  A proof that these multiplicity functions satisfy the consistency equation can be found in \cite{BRS}.

Now, by Merrill's Theorem, we need to find a set $E$ that satisfies the conditions of the theorem.  To construct $E$, we take the support of $m_k^2$, call it $E_{m_k^2}$ and put that in $E$.  Now, we take the interval where $m_k^2(\omega) \geq k - j$ for $j = 1, 2, \ldots k-2$ and shift the negative half, $[ \frac{-2^j\pi}{2^k - 1}, 0)$ to the right $2^j\pi$, and call this $E_j^+$.  We then take the positive half, $[0, \frac{2^j\pi}{2^k -1})$ and shift it to the left $2^j\pi$ and call this part $E_j^-$.  Then let $E_j = E_j^- \cup E_j^+$.

Define, \[E = E_{m_k^2} \cup \bigcup_{j=1}^{k-2} E_j;\]
we have:
\begin{multline*}
E = \cup_{j=1}^{k-2}[-2^j\pi, -2^j\pi + \frac{2^j\pi}{2^k - 1}) \cup [-\pi, \frac{-(2^k - 2)\pi}{2^k - 1}) \cup [\frac{-2^{k-1}\pi}{2^k - 1}, \frac{2^{k-1}}{2^k - 1}\pi) \\
\cup [ \frac{(2^k - 2)\pi}{2^k - 1}, \pi) \cup \cup_{j=1}^{k-2}[2^j\pi - \frac {2^j\pi}{2^k - 1}, 2^j\pi).
\end{multline*}
By definition, $E$ satisfies the first condition of Merrill's theorem.  The other two conditions are easily verified.  Hence, there is a wavelet set $W = 2E \setminus E$:
\begin{multline*}
W = [-2^{k-1}\pi, \frac{-(2^{2k-1}-2^k)\pi}{2^k - 1}) \cup [\frac{-2^k\pi}{2^k -1}, -\pi) \cup [\frac{-(2^k - 2)\pi}{2^k -1}, \frac{-2^{k-1}\pi}{2^k -1}) \\
\cup [\frac{2^{k-1}\pi}{2^k -1}, \frac{(2^k - 2)\pi}{2^k -1}) \cup[\pi, \frac{2^k\pi}{2^k -1}) \cup [\frac{(2^{2k-1} - 2^k)\pi}{2^k -1}, 2^{k-1}\pi).
\end{multline*}

\begin{proof}[Proof of theorem~\ref{T:Main} (for dilation factor 2)]
Now we shall construct a second wavelet set to form an interpolation pair.  Take the same multiplicity function, $m_k^2(\omega)$ and construct an $E'$ similar to $E$, except we will take the interval $[\frac{-2^{k-1}}{2^k -1},\frac{-\pi}{2})$ and shift it right $2^{k-1}\pi$.  The result is
\begin{multline*}
E' = \cup_{j=1}^{k-2}[-2^j\pi, -2^j\pi + \frac{2^j\pi}{2^k - 1}) \cup [-\pi, \frac{-(2^k - 2)\pi}{2^k - 1}) \\
\cup [\frac{-\pi}{2}, \frac{2^{k-1}}{2^k - 1}) \cup [ \frac{(2^k - 2)\pi}{2^k - 1}, \pi) \cup \cup_{j=1}^{k-2}[2^j\pi - \frac{2^j\pi}{2^k - 1}, 2^j\pi) \\
\cup [ 2^{k-1}\pi -\frac{2^{k-1}}{2^k -1}, 2^{k-1}\pi - \frac{\pi}{2}).
\end{multline*}
This new $E'$ also satisfies the conditions of Merrill's theorem, so there is a wavelet set $W'$:
\begin{multline*}
W' = [-2^{k-1}\pi, \frac{-(2^{2k-1}-2^k)\pi}{2^k - 1}) \cup [\frac{-(2^k - 2)\pi}{2^k - 1}, \frac{-\pi}{2}) \\
\cup [\frac{2^{k-1}\pi}{2^k -1}, \frac{(2^k - 2)\pi}{2^k -1}) \cup [\pi, \frac{2^k\pi}{2^k -1}) \cup [2^{k-1}\pi - \frac{\pi}{2}, 2^{k-1}\pi) \\
\cup [2^k\pi - \frac{2^k\pi}{2^k -1}, 2^k\pi - \pi)
\end{multline*}
\begin{multline*}
W \cap W' = [-2^{k - 1}\pi, \frac{-(2^{2k}-2^k)\pi}{2^k - 1}) \cup [\frac{-(2^k - 2)\pi}{2^k - 1}, \frac{-2^{k-1}\pi}{2^k - 1}) \\
\cup [\frac{2^{k-1}\pi}{2^k - 1}, \frac{(2^k - 2)\pi}{2^k - 1}) \cup [ \pi, \frac{2^k \pi}{2^k - 1}) \cup [2^{k-1}\pi - \frac{\pi}{2}, 2^{k-1}\pi)
\end{multline*}
\[ W \setminus W' = [\frac{-2^k \pi}{2^k - 1}, - \pi) \cup [\frac{(2^{2k-1} - 2^k)\pi}{2^k - 1}, 2^{k-1}\pi - \frac{\pi}{2}) = A_1 \cup A_2 \]
\[ W' \setminus W = [\frac{- 2^{k - 1} \pi}{2^k - 1}, \frac{-\pi}{2}) \cup [2^k \pi - \frac{2^k \pi}{2^k - 1}, 2^k \pi - \pi) = B_1 \cup B_2 \]
We claim that the wavelet sets $W$ and $W'$ form an interpolation pair.  Then $\sigma: W \to W'$ is defined as follows:
\begin{equation*}
\sigma(\xi) = 
\begin{cases}
\xi, & \xi \in W \cap W'\\
\xi + 2^k\pi, & \xi \in [\frac{-2^k\pi}{2^k - 1}, -\pi), \hfill (A_1 \to B_2) \\
\xi - 2^{k-1}\pi, & \xi \in [\frac{(2^{2k-1} - 2^k)\pi}{2^k -1}, 2^{k-1} -\frac{\pi}{2}), \hfill (A_2 \to B_1)
\end{cases}
\end{equation*}
We need to check that $\sigma$ is involutive (see \cite{DL}).  Notice that $B_2 = 2A_2$ and that $2B_1 = A_1$.  For $\xi \in A_1, \: \sigma(\xi) \in B_2, \text{ so } \frac{1}{2} \sigma(\xi) \in A_2$.  This yields the following:
\begin{align*}
\sigma^2(\xi) & = 2 \sigma(\frac{1}{2} \sigma(\xi)) \\
& = 2 \sigma(\frac{1}{2}(\xi + 2^k \pi)) \\
& = 2 \sigma(\frac{1}{2}\xi + 2^{k-1} \pi) \\
& = 2(\frac{1}{2} \xi + 2^{k-1} \pi - 2^{k -1} \pi) \\
& = \xi
\end{align*}
Similarly, for $\xi \in A_2$, we get that $\sigma^2(\xi) = \xi$, whence $\sigma$ is involutive.
Finally, to interpolate between these two wavelets, define the functions $h_1$ on $W$ and $h_2$ on $W'$ as follows:
\begin{align*}
h_1(\xi) &=
\begin{cases}
1, & \xi \in W \cap W', \\
\frac{1}{\sqrt{2\pi}}, & \xi \in [\frac{(2^{2k-1} - 2^k)\pi}{2^k -1}, 2^{k-1} -\frac{\pi}{2}), \\
\frac{1}{\sqrt{2\pi}}, & \xi \in [\frac{-2^k\pi}{2^k - 1}, -\pi),
\end{cases}
\\
h_2(\xi) &=
\begin{cases}
1, & \xi \in W \cap W', \\
\frac{1}{\sqrt{2\pi}}, & \xi \in [-\frac{(2^k-2^{k-1})\pi}{2^k -1}, -\frac{\pi}{2}), \\
-\frac{1}{\sqrt{2\pi}}, & \xi \in [2^k\pi-\frac{2^k\pi}{2^k-1},2^k\pi-\pi),
\end{cases}
\end{align*}
and extend via 2-dilation periodicity.  It is also routine to verify that these functions satisfy the matrix condition~\ref{E:M}.

\begin{claim}
Let $W$, $W'$, $h_1$ and $h_2$ be as above.  Then the wavelet defined by
\[\hat{\psi} = h_1 \hat{\psi}_{W} + h_2 \hat{\psi}_{W'}\]
is an element of $\Cal{M}_{k-2}$
\end{claim}
\begin{proof}
As defined above, $\hat{\psi}$ is a wavelet.  Note that $\sigma$ above yields that $2^{k-2},2^{k-1} \in \Cal{E}_{\psi}$, whence the claim is established by corollary~\ref{C:1}.
\end{proof}
Consequently, there exist wavelets in all $\Cal{M}_n$'s, for dilation by 2.
\end{proof}

\begin{example} [k=3] \label{E:1}
The following is a concrete example of the above procedure.  The resulting wavelet is in $\Cal{M}_1$, so that the corresponding subspace $V_{-1}$ is invariant under translations.  Additionally, this example will be of use in section~\ref{S:T}.

The multiplicity function is:
\begin{equation*}
m(\omega) = 
\begin{cases}
2, & \omega \in [-\frac{2\pi}{7}, \frac{2\pi}{7}), \\
0, & \omega \in [-\frac{6\pi}{7}, -\frac{4\pi}{7}) \cup [\frac{4\pi}{7}, \frac{6\pi}{7}), \\
1, & \text{for all other $\omega \in [-\pi,\pi)$}.
\end{cases}
\end{equation*}
Note that this is the multiplicity function for the Journ\'{e} wavelet.
\begin{align*}
W &= [-4\pi, \frac{-24\pi}{7}) \cup [\frac{-8\pi}{7}, -\pi) \cup [\frac{-6\pi}{7}, \frac{-4\pi}{7}) \\
& \cup [\frac{4\pi}{7}, \frac{6\pi}{7}) \cup [\pi, \frac{8\pi}{7}) \cup [\frac{24\pi}{7}, 4\pi)
\end{align*}
\begin{align*}
W' &= [-4\pi, \frac{-24\pi}{7}) \cup [\frac{-6\pi}{7}, \frac{-\pi}{2}) \cup [\frac{4\pi}{7}, \frac{6\pi}{7}) \\
& \cup [\pi, \frac{8\pi}{7}) \cup [\frac{7\pi}{2}, 4\pi) \cup [\frac{48\pi}{7}, 7\pi)
\end{align*}
We have:
\begin{equation*}
\sigma(\xi) = 
\begin{cases}
\xi, &  \xi \in W \cap W', \\
\xi - 4\pi, & \xi \in [\frac{24\pi}{7}, \frac{7\pi}{2}), \\
\xi + 8\pi, & \xi \in [\frac{-8\pi}{7}, -\pi),
\end{cases}
\end{equation*}
therefore interpolation as above yields a wavelet in $\Cal{M}_1$.
\end{example}

We have found wavelet sets for dilation by 2 and interpolated between them.  Now we would like to generalize this argument 
and be able to do this for any dilation.  It turns out that 2 is the only dilation factor that has symmetric multiplicity functions, the multiplicity functions associated to the others factors are asymmetric.
We will now look at integer dilation factors greater than or equal to 3, such factors will be denoted by $d$.

\begin{theorem} There exists a family of multiplicity functions on $[-\pi, \pi)$, $\{ m_k^{d} \}$ for $ k \geq 2$, $ d \geq 3$ that satisfy
the consisten
cy equation~\ref{Eq:M}, given as follows:
\begin{equation*}
m_k^{d}(\omega) = 
\begin{cases}
0, &  - \pi \leq \omega < \frac{-d^{k - 2}(d- 1)(2\pi)}{d^k - 1}, \hfill (I_1)\\
k - j, & \frac{-d^{j - 1}(d - 1)(2\pi)}{d^k - 1} \leq \omega < \frac{- d^{j - 2}(d - 1)(2 \pi)}{d^k - 1}, \hfill (I_2)\\
k - 1, & \frac{-(d - 1)(2 \pi)}{d^k - 1} \leq \omega < \frac{2 \pi}{d^k - 1}, \hfill (I_3) \\
k - j,  &  \frac{d^{j - 2}(2 \pi)}{d^k - 1} \leq \omega < \frac{d^{j - 1}(2 \pi)}{d^k - 1}, \hfill (I_4)\\
0, & \frac{d^{k - 2}(2\pi)}{d^k - 1} \leq \omega < \frac{2 \pi}{d} - \frac{(d  - 1)(2\pi)}{d(d^k - 1)}, \hfill (I_5)\\
1, & \frac{2\pi}{d} - \frac{(d - 1)(2\pi)}{d(d^k - 1)} \leq \omega < \frac{2\pi}{d} + \frac{2\pi}{d(d^k - 1)}, \hfill (I_6)\\
0, & \frac{2\pi}{d} + \frac{2\pi}{d(d^k - 1)} \leq \omega < \pi, \hfill (I_7)
\end{cases}
\end{equation*}
for, $j = 2, 3, \ldots, k - 1$.  Furthermore, for each of the functions $m_k^{d}$, there is an associated generalized scaling set $E_k^{d}$ and a wavelet set $W_k^{d}$.
\end{theorem}

\begin{proof}   Fix $d$ and $k$.  We shall first give a brief outline to show  that $m_k^{d}(\omega)$ satisfies the consistency equation~\ref{Eq:M}.  For a complete proof of this fact see \cite{Sch}.

Let $\omega \in I_1$, then $\frac{\omega}{d} \in I_2$ where $m(\omega) = k - (k-1)$
and $\frac{\omega}{d} + \frac{2\pi}{d} \in I_5$, and $\frac{\omega}{d} + \frac{4\pi}{d}, \ldots,
\frac{\omega}{d} + \frac{(d - 1)2\pi}{d}$ are in either $I_1$ or $I_7$.  Consistency equation: $0 + 1 = k - (k - 1) + 0 \cdots + 0$.

A similar argument is used for $\omega \in I_m, \: 2 \leq m \leq 7$ to prove that $m_k^d$ satisfies the consistency equation.  Here is how to construct the set $E_k^{d}$.  Let $E_{m_k^{d}}$ = support of $m_k^{d}$ so

\begin{eqnarray*}E_{m_k^{d}} & = & [\frac{-d^{k-2}(d - 1)(2\pi)}{d^k - 1}, \frac{d^{k-2}(2\pi)}{d^k - 1}) \cup \\ & & \mbox{} [\frac{2\pi}{d} - \frac{(d - 1)(2\pi)}{d(d^k - 1)}, \frac{2\pi}{d} + \frac{2\pi}{d(d^k - 1)})
\end{eqnarray*}
Now take the set where $m_k^{d} \geq k - j, \: \text{ for } j = 1, \ldots k - 2$ and shift it to the right by $2\pi(d^{j - 1})$, call this set $E_{j-1}$.  Thus,
\[E_{j - 1}  =  [2\pi(d^{j-1}) - \frac{d^{j-1}(d -1)(2\pi)}{d^k -1},  2\pi(d^{j-1}) + \frac{d^{j-1}(2\pi)}{d^k - 1})\]
and let
\[E_k^{d} = E_{m_k^{d}} \cup \cup_{j=1}^{k-2}E_{j-1}.\]
It is easy to show that \[ E_k^{d} \subseteq d(E_k^{d})\]
so by Merrill's theorem,
\begin{eqnarray*}
W_k^{d} & = & d(E_k^{d}) \setminus E_k^{d} \\
& = & [\frac{-d^{k-1}(d - 1)2\pi}{d^k - 1}, \frac{-d^{k-2}(d - 1)2\pi}{d^k - 1})\\
             &   &  \hspace{.3in} \cup [\frac{d^{k-2}(2\pi)}{d^k - 1}, \frac{2\pi}{d} - \frac{(d - 1)2\pi}{d(d^k - 1)})\\
             &   &  \hspace{.3in} \cup [2\pi(d^{k-2}) - \frac{d^{k-2}(d - 1)2\pi}{d^k - 1}, 2\pi(d^{k-2}) + \frac{d^{k-2}(2\pi)}{d^k - 1})
\end{eqnarray*}
 is the desired wavelet set.
\end{proof}

\begin{proof}[Proof of theorem~\ref{T:Main} (for dilation factor greater than 2)]
The above is one way of constructing $E_k^{d}$ and $W_k^{d}$.  Now, we will construct a second generalized scaling set $E$ and a second wavelet set $W$ that have the same multiplicity function by moving the interval $[\frac{-d^{k-2}(d - 1)2\pi}{d^k - 1}, \frac{-2\pi}{d^2})$ to the right
$2\pi(d^{k-2})$.
\begin{multline*}
\widetilde{E_k^{d}} = [\frac{-2\pi}{d^2}, \frac{d^{k-2}(2\pi)}{d^k - 1}) \cup [\frac{2\pi}{d} - \frac{(d - 1)2\pi}{d(d^k - 1)}, \frac{2\pi}{d} + \frac{2\pi}{d(d^k - 1)}) \\
\cup_{j=1}^{k-2} E_{j-1} \cup [2\pi(d^{k-2}) - \frac{d^{k-2}(d - 1)2\pi}{d^k - 1}, 2\pi(d^{k-2}) - \frac{2\pi}{d^2})
\end{multline*}

\begin{multline*}
\widetilde{W_k^{d}} = [\frac{-2\pi}{d}, \frac{-2\pi}{d^2}) \cup [\frac{d^{k-2}(2\pi)}{d^k - 1}, \frac{2\pi}{d} - \frac{(d - 1)2\pi}{d(d^k - 1)})\\
\cup [2\pi(d^{k-2}) - \frac{2\pi}{d^2}, 2\pi(d^{k-2}) + \frac{d^{k-2}(2\pi)}{d^k - 1})\\
\cup [2\pi(d^{k-1}) - \frac{d^{k-1}(d - 1)2\pi}{d^k - 1}, 2\pi(d^{k-1}) - \frac{2\pi}{d})
\end{multline*}

Consider the wavelet sets $W_k^{d}$ and $\widetilde{W_k^{d}}$, then we have the following:
\begin{multline*}
W_k^{d} \cap \widetilde{W_k^{d}} = [\frac{-2\pi}{d}, \frac{-d^{k-2}(d - 1)2\pi}{d^k - 1}) \\
\cup [\frac{d^{k-2}(2\pi)}{d^k - 1}, \frac{2\pi}{d} - \frac{(d - 1)2\pi}{d(d^k - 1)}) \cup [2\pi d^{k-2} - \frac{2\pi}{d^2}, 2\pi d^{k-2} + \frac{d^{k-2} 2\pi}{d^k - 1}) 
\end{multline*}
\begin{multline*}
W_k^{d} \setminus \widetilde{W_k^{d}} = [\frac{-d^{k-1}(d - 1)2\pi}{d^k - 1}, \frac{-2\pi}{d}) \\
\cup [2\pi d^{k-2} - \frac{d^{k-2}(d - 1)2\pi}{d^k - 1}, 2\pi d^{k-2} - \frac{2\pi}{d^2})
\end{multline*}
\begin{multline*}
\widetilde{W_k^{d}} \setminus W_k^{d} = [\frac{-d^{k-2}(d - 1)2\pi}{d^k - 1}, \frac{-2\pi}{d^2}) \\
\cup [2\pi d^{k-1} - \frac{d^{k-1}(d - 1)2\pi}{d^k - 1}, 2\pi d^{k-1} - \frac{2\pi}{d})
\end{multline*}

This gives us the following map $\sigma(\xi): W_k^{d} \rightarrow \widetilde{W_k^{d}}$:
\begin{equation*}
\sigma(\xi) = 
\begin{cases}
\xi, & \xi \in W_k^d \cap \widetilde{W_k^d}, \\
\xi + d^{k-1}2\pi, & \xi \in [\frac{-d^{k-1}(d - 1)2\pi}{d^k - 1}, \frac{-2\pi}{d}), \\
\xi - d^{k-2}2\pi, & \xi \in [2\pi d^{k-2} - \frac{d^{k-2}(d - 1)2\pi}{d^k - 1}, 2\pi d^{k-2} - \frac{2\pi}{d^2}).
\end{cases}
\end{equation*}
As with the case of dilation by two, it is routine to check that $\sigma$ is involutive.  Additionally, we need to define functions $h_1$ and $h_2$ for the interpolation;  define the functions $h_1$ on $W_k^{d}$ and $h_2$ on $\widetilde{W_k^{d}}$ as follows:
\begin{align*}
h_1(\xi) &=
\begin{cases}
1, & \xi \in W_k^{d} \cap \widetilde{W_k^{d}}, \\
\frac{1}{\sqrt{2\pi}}, & \xi \in [\frac{-d^{k-1}(d - 1)2\pi}{d^k - 1}, \frac{-2\pi}{d}), \\
\frac{1}{\sqrt{2\pi}}, & \xi \in [2\pi d^{k-2} - \frac{d^{k-2}(d - 1)2\pi}{d^k - 1}, 2\pi d^{k-2} - \frac{2\pi}{d^2}),
\end{cases}
\\
h_2(\xi) &=
\begin{cases}
1, & \xi \in W_k^{d} \cap \widetilde{W_k^{d}}, \\
\frac{1}{\sqrt{2\pi}}, & \xi \in [\frac{-d^{k-2}(d - 1)2\pi}{d^k - 1}, \frac{-2\pi}{d^2}), \\
-\frac{1}{\sqrt{2\pi}}, & \xi \in [2\pi d^{k-1} - \frac{d^{k-1}(d - 1)2\pi}{d^k - 1}, 2\pi d^{k-1} - \frac{2\pi}{d}),
\end{cases}
\end{align*}
and extend via d-dilation periodicity.  

\begin{claim}
Let $W_k^{d}$ and $\widetilde{W_k^{d}}$, $h_1$ and $h_2$ be as above.  Then the wavelet defined by
\[\hat{\psi} = h_1 \hat{\psi}_{W} + h_2 \hat{\psi}_{W'}\]
is an element of $\Cal{M}_{k-2}$
\end{claim}
\begin{proof}
As defined above, $\hat{\psi}$ is a wavelet.  Note that $\sigma$ above yields that $d^{k-2},d^{k-1} \in \Cal{E}_{\psi}$, whence the claim is established by corollary~\ref{C:1}.
\end{proof}
Hence, there exist wavelets in all $\Cal{M}_n$'s, for dilation by $d > 2$.  We have now completed the proof of theorem~\ref{T:Main}.
\end{proof}

\begin{example}: $d = 3, k = 4$
\begin{equation*}
m_4^{3}(\omega) = 
\begin{cases}
0, &  - \pi \leq \omega < \frac{-18\pi}{40}, \\
1, & \frac{-18\pi}{40} \leq \omega < \frac{-6\pi}{40}, \\
2, & \frac{-6\pi}{40} \leq \omega < \frac{-2 \pi}{40}, \\
3, & \frac{-2 \pi}{40} \leq \omega < \frac{\pi}{40}, \\
2, & \frac{\pi}{40} \leq \omega < \frac{3 \pi}{40}, \\
1,   &  \frac{3\pi}{40} \leq \omega < \frac{9 \pi}{40}, \\
0, & \frac{9\pi}{40} \leq \omega < \frac{26 \pi}{40}, \\
1, & \frac{26\pi}{40} \leq \omega < \frac{27\pi}{40}, \\
0, & \frac{27\pi}{40} \leq \omega < \pi.
\end{cases}
\end{equation*}

From the construction in the proof above, we get:
\[
E_4^3 = [\frac{-18\pi}{40}, \frac{9\pi}{40}) \cup [\frac{26\pi}{40}, \frac{27\pi}{40})
\cup [\frac{78\pi}{40},\frac{81\pi}{40}) \cup [\frac{234\pi}{40}, \frac{243\pi}{40})
\]
Since
\[W_4^3 = 3E_4^3 \setminus E_4^3, \]
the resulting wavelet set is:
\[W_4^3 = [\frac{-54\pi}{40}, \frac{-18\pi}{40}) \cup [\frac{9\pi}{40}, \frac{26\pi}{40}) \cup [\frac{702\pi}{40}, \frac{729\pi}{40}).\]
\par
If we take the same multiplicity function and do the second construction, we get:
\begin{multline*}
\widetilde{E_4^3} = [\frac{-2\pi}{9}, \frac{9\pi}{40}) \cup [\frac{26\pi}{40}, \frac{27\pi}{40}) \cup [\frac{78\pi}{40}, \frac{81\pi}{40}) \\
\cup [\frac{234\pi}{40}, \frac{243\pi}{40}) \cup [\frac{702\pi}{40}, \frac{160\pi}{9})
\end{multline*}
The second wavelet set is:
\[
\widetilde{W_4^3} = [\frac{-2\pi}{3}, \frac{-2\pi}{9}) \cup [\frac{9\pi}{40}, \frac{26\pi}{40}) \cup [\frac{160\pi}{9},\frac{729\pi}{40})
\cup [\frac{2106\pi}{40}, \frac{160\pi}{3}).
\]

Now we want to find a map $\sigma(\xi): W_4^3 \rightarrow \widetilde{W_4^3}$.  First, we will
identify the intersection of the wavelet sets.
\begin{gather*}
W_4^3 \cap \widetilde{W_4^3} = [\frac{-2\pi}{3}, \frac{-18\pi}{40}) \cup [\frac{9\pi}{40}, \frac{26\pi}{40}) \cup [\frac{160\pi}{9}, \frac{729\pi}{40}) \\
W_4^3 \setminus \widetilde{W_4^3} = [\frac{-54\pi}{40}, \frac{-2\pi}{3}) \cup [\frac{702\pi}{40}, \frac{160\pi}{9}) \\
\widetilde{W_4^3} \setminus W_4^3 = [\frac{-18\pi}{40}, \frac{-2\pi}{9}) \cup [\frac{2106\pi}{40}, \frac{160\pi}{3})
\end{gather*}
We have:
\begin{equation*}
\sigma(\xi) = 
\begin{cases}
\xi, &  \xi \in W_4^3 \cap \widetilde{W_4^3},\\
\xi + 54\pi, & \xi \in [\frac{-54\pi}{40}, \frac{-2\pi}{3}),\\
\xi - 18\pi, & \xi \in [\frac{702\pi}{40}, \frac{160\pi}{9}),
\end{cases}
\end{equation*}
Hence, by interpolating as above, we get a wavelet that is an element of $\Cal{M}_2$ for dilation by 3.
\end{example}

\section{Topology and Concluding Remarks} \label{S:T}
We have shown the existence of wavelets with higher order translation invariance properties.  Not much is known about them, other than what is known about the MSF wavelets ($\Cal{M}_{\infty}$).  For instance, the MSF wavelets are relatively closed in the collection of all wavelets in the $\ltwo$ norm.  We have the following analogue of the this fact: the collections $\Cal{L}_n$ are relatively closed, but the collections $\Cal{M}_n$ are not.  We shall conclude the paper with some possible research directions.

\begin{theorem}
Suppose that $\psi_n$ is a sequence of wavelets that converges to a wavelet $\psi$.  Suppose further that $\{\psi_n\}_n \subset \Cal{L}_j$.  Then $\psi \in \Cal{L}_j$.
\end{theorem}
\begin{proof}
We shall prove by contrapositive.  Suppose that $\psi \notin \Cal{L}_j$, then there exists an $l \in \Cal{E}_{\psi}$ that is not divisible by $d^j$.  We shall show that there exists a positive integer $N$ such that for $n \geq N$, $\psi_n \notin \Cal{L}_j$, by showing that $l \in \Cal{E}_{\psi_n}$.  Note that $\psi_n$ converges to $\psi$ in measure.

Let $F \subset supp(\hat{\psi})$ be such that $F+2\pi l \subset supp(\hat{\psi})$.  Choose $G \subset F$, a set of non-zero measure, such that for $\xi \in G \cup G+2\pi l$, $|\hat{\psi}(\xi)| > \epsilon_1$ for some $\epsilon_1 > 0$.  Let $\epsilon_2$ be the measure of $G$.  Choose $\epsilon < min(\frac{\epsilon_1}{2},\frac{\epsilon_2}{2})$.  Let $N$ be such that for $n \geq N$, $m(\{\xi:|\hat{\psi}(\xi)-\hat{\psi}_n(\xi)| \geq \epsilon\}) < \epsilon$.  Let $H = \{\xi \in G:|\hat{\psi}(\xi)-\hat{\psi}_n(\xi)| < \epsilon\}$; note that $m(H) > m(G) - \epsilon > \frac{m(G)}{2}$.  Likewise, define $H' = \{\xi \in G+2\pi l:|\hat{\psi}(\xi)-\hat{\psi}_n(\xi)| < \epsilon\}$; $m(H') > \frac{m(G)}{2}$.

We have that $H+2\pi l \subset G+2\pi l$ and $H' \subset G+2\pi l$, and both sets have measure greater than half of $G+2\pi l$, whence they must intersect.  Hence, if $H$ and $H'$ are subsets of $supp(\hat{\psi}_n)$, then $\hat{\psi}_n \notin \Cal{L}_j$.

For $\xi \in H \cup H'$, 
\[|\hat{\psi}(\xi)|-|\hat{\psi}_n(\xi)| \leq |\hat{\psi}(\xi)-\hat{\psi}_n(\xi)|< \epsilon < \epsilon_1 < |\hat{\psi}(\xi)|,\]
therefore, $|\hat{\psi}_n(\xi)| > 0$, as required.
\end{proof}

We shall now demonstrate that the $\Cal{M}_n$'s are not closed.  In the following example, we show a sequence of wavelets in $\Cal{M}_1$ that converges to a MSF wavelet in $\ltwo$.  This technique can be extended to any of the interpolation pairs we used in this paper.

Recall from Example~\ref{E:1} that a wavelet was constructed out of an interpolation pair and the resulting wavelet was in $\Cal{M}_1$.  By perturbing one of the wavelet sets in an appropriate manner, a sequence of interpolation pairs can be constructed.  The end result is a sequence of wavelets in $\Cal{M}_1$ which converges to an MSF wavelet.

We shall use the first wavelet set unchanged:
\begin{multline*}
W = [-4\pi, \frac{-24\pi}{7}) \cup [\frac{-8\pi}{7}, -\pi) \cup [\frac{-6\pi}{7}, \frac{-4\pi}{7}) \\
\cup [\frac{4\pi}{7}, \frac{6\pi}{7}) \cup [\pi, \frac{8\pi}{7}) \cup [\frac{24\pi}{7}, 4\pi).
\end{multline*}
Let $a_n$ be a sequence such that $-\frac{4\pi}{7}<a_n<-\frac{\pi}{2}$ and $\lim_{n\to\infty} a_n = -\frac{4\pi}{7}$.  Consider the sequence of generalized scaling sets (they satisfy the conditions of Merrill's theorem):
\begin{multline*}
E_n = [-2\pi, \frac{-12\pi}{7}) \cup [-\pi, \frac{-6\pi}{7}) \cup [a_n, \frac{4\pi}{7}) \\
\cup [\frac{6\pi}{7}, \pi) \cup [\frac{12\pi}{7}, 2\pi) \cup [\frac{24\pi}{7}, 4\pi + a_n).
\end{multline*}
The corresponding wavelet sets are:
\begin{multline*}
W_n = [-4\pi, \frac{-24\pi}{7}) \cup [2a_n, -\pi) \cup [\frac{-6\pi}{7}, a_n) \\
\cup [\frac{4\pi}{7}, \frac{6\pi}{7}) \cup [\pi, \frac{8\pi}{7}) \cup [4\pi +a_n, 4\pi) \cup [\frac{48\pi}{7}, 8\pi + 2a_n).
\end{multline*}
It is easy to verify that for all $n$, $W$ and $W_n$ form an interpolation pair, where $\sigma$ is given by:
\begin{equation*}
\sigma(\xi) = 
\begin{cases}
\xi, &  \xi \in W \cap W_n\\
\xi - 4\pi, & \xi \in [\frac{24\pi}{7}, 4\pi + a_n) \\
\xi + 8\pi, & \xi \in [\frac{-8\pi}{7}, 2a_n)
\end{cases}
\end{equation*}
as in example~\ref{E:1}.  Note that $W_n \to W$ as $n \to \infty$, in the symmetric difference metric, whence the limit of the sequence of interpolated wavelets is precisely $\psi_W$.

There are many open questions regarding wavelets with higher order translation invariance properties.  The first is whether there exist ``nice'' wavelets with these properties.  The wavelets constructed in this paper are non-MRA wavelets and are discontinuous in the frequency domain.  As such, these particular wavelets are not suitable for applications.  A second consideration is the extension of this to several dimensions.  The characterization theorem (\ref{C:1}) has a natural extension, but operator interpolation in several dimensions is difficult as finding interpolation pairs is non-trivial.

\bibliographystyle{amsplain}
\bibliography{wtipon}
\nocite{HW}
\nocite{P}
\nocite{FW}
\nocite{C}
\nocite{W1}
\end{document}